\documentclass[reqno,centertags, 12pt]{amsart}
\usepackage{amsmath,amsthm,amscd,amssymb}
\usepackage{latexsym}
\newcommand{\bbR}{{\mathbb{R}}}

\newcommand{\bbZ}{{\mathbb{Z}}}
\newcommand{\bbC}{{\mathbb{C}}}

\newcommand{\dott}{\,\cdot\,}

\newcommand{\lb}{\label}
\newcommand{\f}{\frac}

\newcommand{\ol}{\overline}

\newcommand{\dist}{\text{\rm{dist}}}

\newcommand{\spec}{\text{\rm{spec}}}

\newcommand{\supp}{\text{\rm{supp}}}

\newcommand{\bi}{\bibitem}

\newcommand{\beq}{\begin{equation}}
\newcommand{\eeq}{\end{equation}}
\newcommand{\ba}{\begin{align}}
\newcommand{\ea}{\end{align}}
\newcommand{\veps}{\varepsilon}
\newcounter{smalllist}
\newenvironment{SL}{\begin{list}{{\rm\roman{smalllist})}}{%
\setlength{\topsep}{0mm}\setlength{\parsep}{0mm}\setlength{\itemsep}{0mm}%
\setlength{\labelwidth}{2em}\setlength{\leftmargin}{2em}\usecounter{smalllist}%
}}{\end{list}}

\allowdisplaybreaks
\numberwithin{equation}{section}
\newtheorem{theorem}{Theorem}[section]
\newtheorem*{t0}{Theorem}
\newtheorem*{t1}{Theorem 1}
\newtheorem*{t2}{Theorem 2}

\newtheorem*{p2.1}{Proposition 2.1}
\newtheorem{proposition}[theorem]{Proposition}
\newtheorem{lemma}[theorem]{Lemma}

\theoremstyle{definition}

\theoremstyle{remark}

\newcommand{\abs}[1]{\lvert#1\rvert}

\begin{document}
\title[Zeros of Orthogonal Polynomials on the Real Line]{Zeros of Orthogonal Polynomials \\
on the Real Line}
\author[S. Denisov and B. Simon]{Sergey A.~Denisov$^1$ and Barry Simon$^{1,2}$}

\thanks{$^1$ Mathematics 253-37, California Institute of Technology, Pasadena, CA 91125. 
E-mail: denissov@its.caltech.edu; bsimon@caltech.edu}
\thanks{$^2$ Supported in part by NSF grants DMS-9707661, DMS-0140592}

\date{July 5, 2002}

\begin{abstract} Let $p_n(x)$ be the orthonormal polynomials associated to a measure $d\mu$ of compact 
support in $\bbR$. If $E\notin\supp(d\mu)$, we show there is a $\delta>0$ so that for all $n$, either 
$p_n$ or $p_{n+1}$ has no zeros in $(E-\delta, E +\delta)$. If  $E$ is an isolated point of 
$\supp(\mu)$, we show there is a $\delta$ so that for all $n$, either $p_n$ or $p_{n+1}$ has at 
most one zero in $(E-\delta, E+\delta)$. We provide an example where the zeros of $p_n$ are dense in 
a gap of $\supp (d\mu)$. 
\end{abstract}

\maketitle

\section{Introduction} \lb{s1}

Let $d\mu$ be a measure on $\bbR$ whose support is not a finite number of points and with $\int 
\abs{x}^n \, d\mu (x) < \infty$ for all $n=0,1,2,\dots$. The orthonormal polynomials $p_n(x;d\mu)$ 
or $p_n(x)$ are determined uniquely by  
\begin{gather}
p_n(x) = \gamma_n x^n + \text{lower order}  \qquad \gamma_n >0 \lb{1.1} \\ 
\int p_n(x) p_m (x)\, d\mu(x) =\delta_{nm} \lb{1.2} 
\end{gather}
There are $a_n >0$, $b_n\in\bbR$ for $n\geq 1$ so that 
\begin{equation} \lb{1.3}
xp_n(x) = a_{n+1} p_{n+1}(x) + b_{n+1} p_n(x) + a_n p_n (x)
\end{equation}
(many works use $a_{n-1}$, $b_{n-1}$ where we use $a_n$, $b_n$). 

In this paper, we will be interested in the zeros of $p_n(x;d\mu)$. The following results are classical 
(see, e.g., Freud's book \cite{FrBk}):
\begin{SL}
\item[(1)] The zeros of $p_n(x)$ are real and simple. 
\item[(2)] If $(a,b)\cap\supp (d\mu)=\emptyset$, then if $a=-\infty$ or $b=+\infty$, $p_n$ has no 
zeros in $(a,b)$ and, in any event, $(a,b)$ has at most one zero of $p_n(x)$. 
\item[(3)] In the determinate case, if $x_0\in\supp(d\mu)$ and $\delta >0$, for all large $n$, $p_n(x)$ 
has a zero in $(x_0-\delta, x_0 + \delta)$. 
\end{SL}

Define 
\[
N_n (x_0, \delta) = \#\text{ of zeros of $p_n(x)$ in } (x_0 -\delta, x_0 +\delta) 
\]
Then (1)--(3) immediately imply:
\begin{SL}
\item[(i)] If $x_0$ is a non-isolated point of $\supp (d\mu)$, then for any $\delta >0$, 
$\lim_{n\to\infty} N_n (x_0,\delta) =\infty$. 
\item[(ii)] If $x_0$ is an isolated point of $\supp (d\mu)$ and $\delta =\dist (x_0, \supp(d\mu) 
\backslash\{x_0\})$, then $N_n (x_0, \delta)$ is never more than $2$, and for all $\delta >0$ and $n$ 
large, $N_n (x_0,\delta)\geq 1$. 
\item[(iii)] If $x_0\notin\supp (d\mu)$ and $\delta=\dist (x_0, \supp (d\mu))$, then $N_n (x_0, 
\delta)$ is never more than $1$. 
\end{SL}

(i) is fairly complete, but (ii), (iii) leave open how often there is one vs.~two points in case 
(ii) and zero vs.~one in case (iii). One might guess that a zero near $x_0\notin\supp(d\mu)$ and 
two zeros near an isolated $x_0$ in $\supp(d\mu)$ are not too common occurrences. 

\smallskip
\noindent{\bf Example.} If $d\mu$ is even about $x=0$, then $p_n(-x)=(-1)^n p_n(x)$. Thus, if 
$n$ is odd, $p_n(0)=0$. So if $0\notin\supp (d\mu)$, we still have $N_n (0,\delta)=1$ for 
all small $\delta$ and $n$ odd. If zero is an isolated point of $d\mu$, $p_n$ for $n$ even has 
a zero at $x_n$ near $0$, but not equal to $0$ (since zeros are simple), so also at $-x_n$, that is, 
$N_n (0,\delta)=2$ for $\delta$ small and $n$ even. So ``not too common" can be as often as $50\%$ 
of the time. Our goal here is to show this $50\%$ is a maximal value. 

\smallskip
It is surprising that there do not seem to be any results on these issues until a recent paper of 
Ambroladze \cite{Amb}, who proved 

\begin{t0}[Ambroladze \cite{Amb}] If $\supp (d\mu)$ is bounded and $x_0\notin\supp (d\mu)$, then for some 
$\delta >0$, $\liminf_{n\to\infty} N_n (x_0,\delta) =0$. 
\end{t0} 

Thus we can use $N_n (x_0,\delta)$ to distinguish when $x_0\in\supp (d\mu)$. Our goal in this paper 
is to prove 

\begin{t1} Let $d=\dist (x_0, \supp(d\mu)) >0$. Let $\delta_n = d^2/(d+\sqrt{2}\, a_{n+1})$ 
{\rm{(}}where $a_n$ is the recursion coefficient given by \eqref{1.3}{\rm{)}}. Then either $p_n$ or 
$p_{n+1}$ {\rm{(}}or both{\rm{)}} has no zeros in $(x_0-\delta_n, x_0 + \delta_n)$. In particular, if 
$a_\infty =\sup_n a_n < \infty$ and $d_\infty = d^2 /(d+\sqrt{2}\, a_\infty)$, then $(x_0 -
\delta_\infty, x_0 + \delta_\infty)$ does not have zeros of $p_j$ for two successive values of $j$. 
\end{t1} 

\begin{t2} Let $x_0$ be an isolated point of $\supp (d\mu)$. Then there exists a $d_0 >0$, so that 
if $\delta_n = d_0^2 /(d_0 + \sqrt{2}\, a_{n+1})$, then at least one of $p_n$ and $p_{n+1}$ has no 
zeros  or one zero in $(x_0 -\delta_n, x_0 + \delta_n)$. In particular, if $a_\infty = \sup_n a_n < 
\infty$ and $\delta_\infty =d_0^2/(d_0 + \sqrt2\, a_\infty)$, then for all large $n$, either $N_n 
(x_0, \delta_\infty) =1$ or $N_{n+1}(x_0, \delta_\infty) =1$. 
\end{t2} 

We will prove Theorem~1 in Section~2 and Theorem~2 in Section~3. In Section~4, we present an example 
of a set of polynomials whose zeros are dense in a gap of the spectrum. 

It is a pleasure to thank Leonid Golinskii and Paul Nevai for useful correspondence. 

\bigskip
\section{Points Outside the Support of $d\mu$} \lb{s2} 

We arrived at the following lemma by trying to abstract the essence of Ambroladze's argument 
\cite{Amb}; it holds for orthogonal polynomials on the complex plane. Let $d\mu$ be a measure on 
$\bbC$ with finite moments and infinite support, and let $p_n(z;d\mu)$ be the orthonormal polynomials. 
Define the reproducing kernel
\begin{equation} \lb{2.1}
K_n (z,w) = \sum_{j=0}^n p_n(z)\, \ol{p_n(w)}
\end{equation}
so in $L^2 (\bbC, d\mu)$, for any polynomial $\pi$ of degree $n$ or less, 
\begin{equation} \lb{2.2}
\int K_n (z,w) \pi(w) \, d\mu (w) = \pi(z)
\end{equation}

\begin{lemma} \lb{L2.1} Suppose $z_0\in\bbC$, $p_j(w)=0$ for some $j\leq n+1$. Then 
\begin{equation} \lb{2.3}
\abs{z_0 -w} \geq \f{\abs{p_j (z_0)}}{K_n (z_0, z_0)^{1/2}} \, \dist (w, \supp (d\mu))
\end{equation}
\end{lemma}

\begin{proof} Let $q(z) = p_j(z)/(z-w)$, which has $\deg (q)\leq n$. Thus, by \eqref{2.2}, $\langle 
K(\dott, z_0), q(\dott)\rangle = q(z_0)$ so, by the Schwarz inequality, 
\[
\f{\abs{p_j (z_0)}}{\abs{z_0-w}} \leq \|q\| \, \|K(\dott, z_0)\|
\]
By \eqref{2.2}, $\|K(\dott, z_0)\| = K(z_0, z_0)^{1/2}$ and clearly, $\|q\|\leq \dist (w, \supp
(d\mu))^{-1} \|p_j\|=\dist (w, \supp(d\mu))^{-1}$. This yields \eqref{2.3}. 
\end{proof}

The following only holds in the real case:

\begin{lemma} \lb{L2.2} For any $x\in\bbR$ and $n$, 
\begin{equation} \lb{2.4}
K_n (x,x) \dist (x,\supp (d\mu))^2 \leq a_{n+1}^2 [p_{n+1}^2 (x) + p_n^2(x)]
\end{equation}
\end{lemma}

\begin{proof} The Christoffel-Darboux formula \cite{FrBk} says 
\[
K_n (x,y) = a_{n+1} \biggl[ \f{p_{n+1} (x) p_n(y) - p_{n+1}(y) p_n(x)}{x-y}\biggr]
\]
so since $\langle p_j, p_k\rangle = \delta_{jk}$, 
\begin{equation} \lb{2.5}
\|(x-\dott) K_n (x,\dott)\|^2 = \abs{a_{n+1}}^2 [p_{n+1}^2 (x) + p_n^2 (x) ]
\end{equation}
Clearly, 
\begin{equation} \lb{2.6}
\|(x-\dott) K_n (x,\dott) \|^2 \geq \dist (x,\supp (d\mu))^2 \|K_n (x, \dott)\|^2
\end{equation} 
and, as above, $\|K_n (x, \dott)\|^2 = K_n (x,x)$, which yields \eqref{2.4}. 
\end{proof} 

{\it Remark.} An alternate way of seeing \eqref{2.5} is to let $\psi$ be the trial vector 
$(p_0(x), \dots, p_n (x),0,0,\dots )$ and note that in terms of the standard Jacobi matrix 
$((J-x)\psi)_j =0$ unless $j=n,n+1$, in which case the values are $-a_{n+1} p_{n+1}(x)$ and 
$a_{n+1} p_n(x)$. \eqref{2.6} is then just $\|(J-x)\psi\|\geq \dist (x, \supp(d\mu)) \|\psi\|$. 

\begin{proof}[Proof of Theorem 1] By \eqref{2.4}, we have that 
\begin{equation} \lb{2.7}
K_n (x_0, x_0) \dist (x_0, \supp (d\mu))^2\leq 2a_{n+1}^2 p_{n+1}^2 (x_0)
\end{equation}
and/or
\begin{equation} \lb{2.8}
K_n (x_0, x_0) \dist (x_0, \supp (d\mu))^2 \leq 2a_{n+1}^2 p_n^2 (x_0)
\end{equation}
Suppose \eqref{2.7} holds. Then, by \eqref{2.3}, if $w$ is a zero of $p_{n+1}(x)$ and if 
$d=\dist (x_0, \supp (d\mu))$, 
\begin{align*}
\abs{x_0-w} &\geq \f{1}{\sqrt2}\, \f{1}{a_{n+1}}\, d \, \dist (w,\supp(d\mu)) \\
&\geq \f{1}{\sqrt2}\, \f{1}{a_{n+1}}\, d(d-\abs{w-x_0})
\end{align*}
which leads directly to $\abs{x_0-w} \geq d^2/(d+a_{n+1}\sqrt2)$. 
\end{proof}

{\it Remark.} There is also a Christoffel-Darboux result for polynomials on the unit circle $\partial 
D = \{z\mid\, \abs{z}=1\}$ in $\bbC$. This leads to the following: If $d\mu$ is a measure on $\partial D$ 
and $z_0\in\partial D$ has $d=\dist (z_0, \supp (d\mu)) > 0$, then the circle of radius $d^2/(2+d)$ 
has no zeros of the orthogonal polynomials. L.~Golinskii has pointed out that the theorem of Fej\'er 
\cite{Fej} that the zeros lie in the convex hull of $\supp (d\mu)$ implies there are no zeros in the 
circle of radius $d^2/2$ --- which is a stronger result, so we do not provide the details. 

\bigskip

\section{Isolated Points of the Support of $d\mu$} \lb{s3} 

To prove Theorem~2, we will make use of the second kind polynomials \cite{FrBk,BS} associated to 
$d\mu$ and $\{p_n\}$. This is a second family of polynomials, $q_n$ defined by recursion coefficients, 
$\tilde a_n$, $\tilde b_n$ with 
\begin{equation} \lb{3.1}
\tilde a_n =a_{n+1} \qquad \tilde b_n = b_{n+1}
\end{equation}
They have the following two critical properties: 

\begin{proposition} \lb{P3.1} 
\begin{SL} 
\item[{\rm{(i)}}] The zeros of $p_{n+1}$ and $q_n$ interlace. In particular, between any two zeros 
of $p_{n+1}$ is a zero of $q_n$. 
\item[{\rm{(ii)}}] If $x_0$ is an isolated point of $d\mu$ and $d\nu$ is a suitable measure with respect 
to which the $q$'s are orthogonal, then $x_0\notin\supp (d\nu)$. 
\end{SL}
\end{proposition}

These are well known. (i) follows from the fact that the zeros of $p_{n+1}$ are eigenvalues of the 
matrix 
\[
J_{ij}^{(n+1)} = b_i \delta_{ij} + a_i \delta_{i\, i+1} + a_{i-1} \delta_{i \, i-1}  
\qquad 1\leq i,j\leq n+1
\]
and the zeros of $q_n$ are the eigenvalues of 
\[
\tilde J_{ij}^{(n)} = \tilde b_i \delta_{ij} + \tilde a_i \delta_{i\, i+1} + \tilde a_{i-1} 
\delta_{i\, i-1} \qquad 1\leq i,j\leq n
\]
which is the matrix $J_{ij}^{(n+1)}$ with the top row and left column removed. (ii) follows because 
of the relation that $\nu$ obeys for all $z\in\bbC\backslash\bbR$ \cite{BS}: 
\begin{equation} \lb{3.2}
\int \f{d\nu(x)}{x-z} = a_1^{-2} \biggl[ b_1 -z - \biggl( \int \f{d\mu(x)}{x-z}\biggr)^{-1}\biggr]
\end{equation}
(if the moment problem is indeterminate, this is one possible $\nu$). Isolated points of $d\mu$ 
are poles of $\int d\mu(x)/(x-z)$ so $\int d\nu(x)/(x-z)$ is regular there. 

\begin{proof}[Proof of Theorem 2] Let $d_0 =\dist (x_0,\supp (d\nu)) >0$ by (ii) of 
Proposition~\ref{P3.1}. By Theorem~1 and \eqref{3.1}, either $q_{n-1}$ or $q_n$ has no zeros in 
$(x_0-\delta_n, x_0 + \delta_n)$. By the intertwining result (Proposition~\ref{P3.1}(i)), either 
$p_n$ or $p_{n+1}$ cannot have two zeros in this interval. 
\end{proof}

{\it Remark.} If $b\in\supp (d\mu)$ is such that $\abs{x_0 -b} =\dist (x_0, \supp (d\mu))$ 
and $\int d\mu(y)/\abs{y-b}=\infty$, then $d\nu$ has an isolated point in between $x_0$ and $b$, 
and so $d_0 < \dist (x_0, \supp (d\mu))$.

\bigskip

\section{An Example of Dense Zeros in the Gap} \lb{s4} 

Nevai raised the issue of whether as $n$ varies, the single possible zero of $p_n$ in a gap $(a,b)$ 
of $\supp (d\mu)$ can yield all of $(a,b)$ as limit points, or if the situation of a single (or 
finite number of) limit point as in the example in Section~\ref{s1} is the only possibility. In this 
section, we describe an explicit bounded Jacobi matrix so that $\supp (d\mu)=[-5,-1]\cup [1,5]$ 
but the set $\{x\in (-1,1)\mid p_n(x)=0 \text{ for some } n\}$ is dense in $[-1,1]$. 

Let $\{\beta_j\}_{j=1}^\infty$ be the sequence 
\[
\beta_1, \beta_2, \ldots = 0, -\tfrac12, 0, \tfrac12, -\tfrac34, -\tfrac12, -\tfrac14, 0, 
\tfrac14, \tfrac12, \tfrac34, -\tfrac78, \dots 
\]
which goes through {\it all} dyadic rationals in $(-1,1)$ with denominator $2^k$ successively for 
$k=1,2,3,\dots$ with each $j/2^k$ ``covered" multiple times. Let $L$ be the Jacobi matrix with 
\begin{alignat}{2}
& a_{2n-1}=3, \qquad a_{2n}=1, &&\qquad n=1,2,\dots \lb{4.1} \\
&b_k = \beta_n &&\qquad \text{if } 2n^2 \leq k < 2(n+1)^2 \lb{4.2} \\
&b_1 = \beta_1 \notag 
\end{alignat}

We claim that 
\begin{SL} 
\item[(1)] $\supp (d\mu) = [-5,-1]\cup [1,5]$ 
\item[(2)] There is an $x_n$ with $\abs{x_n-\beta_n}\leq 2 \ 3^{-2n}$ so that 
\begin{equation} \lb{4.3}
P_{2(n+1)^2 -1} (x_n) =0
\end{equation}
\end{SL}
This provides the claimed example. 

\smallskip
{\it Remarks.} 1. By adjusting $a_1$ and $a_2$ (but keeping $a_{2n+1}=a_1$; $a_{2n} =a_2$), we can 
replace $[-5,-1]\cup [1,5]$ by $[-3-\veps, -1]\cup [1,3+\veps]$, but our method seems to require bands 
larger than the size of the gap. 

\smallskip
2. One can replace \eqref{4.2} by $b_k =\beta_n$ for $\ell_n \leq k < \ell_{n+1}$ so long as 
$\ell_{n+1} -\ell_n \to\infty$. 

\smallskip
3. We believe that the measure associated to $L$ is purely singular. It is perhaps true that the 
phenomenon of zeros dense in a gap requires purely singular spectral measure. 

\medskip
To prove the claims, we let $L_0$ be the Jacobi matrix with $a$'s given by \eqref{4.1} but $b_n=0$,  
and $L_\infty$ the period two, doubly infinite matrix on $\bbZ$ which equals $L_0$ when restricted 
to $\bbZ^+$. By the general theory of periodic Schr\"odinger operators \cite{RS4}, the spectrum of 
$L_\infty$ is the two bands where $\abs{\Delta (x)}\leq 2$ where $\Delta$ is the discriminant, that is, 
the trace of the two-step transfer matrix. If $a_1$, $a_2$ are the two values of $a$ (so $a_1=3$, 
$a_2=1$ in our example), a simple calculation shows that 
\[
\Delta (x) = \f{1}{a_1 a_2}\, (x^2-(a_1^2 + a_2^2))
\]
so $\Delta (x) =\pm 2$ occurs at $x=\pm \abs{a_1 \pm a_2}$. Thus 
\begin{equation} \lb{4.4}
\spec (L_\infty) = [-4, -2] \cup [2, 4]
\end{equation}

The orthonormal polynomials $p_n^{(0)}$ for $L_0$ at $x=0$ obey the recursion relation 
\[
p_{2n+2}^{(0)} (0) = -3 p_{2n}^{(0)}(0)
\]
so we have 
\begin{equation} \lb{4.5}
p_{2n+1}^{(0)}=0 \qquad p_{2n}^{(0)} (0) = (-3)^n
\end{equation}

By the general theory of restricting periodic operators to the half-line, $\spec (L_0)$ is $\spec 
(L_\infty)$ plus a possible single eigenvalue in the gap $(-2,2)$. Since there is a symmetry, the 
only possible eigenvalue is at $x=0$, but \eqref{4.5} says that $0$ is not an eigenvalue since 
$\sum_{j=0}^\infty \abs{P_j (0)}^2 =\infty$. Thus $\spec(L_0) =[-4,-2]\cup [2,4]$ also. $L-L_0$ is 
a diagonal matrix, so it is easy to see $\|L-L_0\|=\sup_j \abs{\beta_{j}}=1$. Thus $\spec(L) 
\subset \cup_{x\in [-1,1]} x + \spec (L_0)=[-5, -1]\cup [1,5]$. On the other hand, since 
the $b$'s are equal to $\beta_j$ on arbitrary long runs, a Weyl vector argument shows that 
\[
\spec(L)\supset \ol{\cup_j \beta_j + \spec(L_0)} = [-5, -1] \cup [1,5]
\]
so claim 1 is proven. 

Let $L_{n;F}$ be the $n\times n$ matrix obtained by taking the first $n$ rows and columns of $L$. 
Then the zeros of $p_n(x)$ are precisely the eigenvalues of $L_{n;F}$ (see \cite[Proposition~5.6]{BS}).  
Let $\varphi_j$ be the $j$ component vector with $(P_0^{(0)}(0), P_1^{(0)}(0), \dots, P_{j-1}^{(0)}(0))$. 
Then if $j$ is odd so $P_j^{(0)}(0)=0$, and we have $L_{0;j,F}\varphi_j =0$. Thus, if $j=2(n+1)^2 -1$, 
\begin{equation} \lb{4.6}
[(L_{j;F}-\beta_n) \varphi_j]_k = (b_k -\beta_n) \varphi_{j,k}
\end{equation}
If $2n^2\leq k\leq 2(n+1)^2-1$, the right-hand side is zero and its absolute value is always less than 
$2\abs{\varphi_{j,k}}$. Thus 
\begin{align*} 
\f{\|(L_{j;F} -\beta_n)\varphi_j\|^2}{\|\varphi_j\|^2} 
& \leq \f{4 \sum_{k=0}^{n^2 -1} 3^{2k}}{\sum_{k=0}^{(n+1)^2 -1} 3^{2k}} \\
&\leq 4 \ 3^{-4n}
\end{align*}
by a simple estimate. Thus $L_{j;F}$ has an eigenvalue within $2 \ 3^{-2n}$ of $\beta_n$, proving claim 2. 

This completes the example.

\bigskip

%%%%%%%%%%%%%%%%%%%%%%%%%%%%%

\end{document}